\documentclass[11pt]{article}

\usepackage{fullpage,amssymb,amsmath}
\usepackage{graphicx}
\usepackage[capitalise]{cleveref}

\title{A Geometric Interpretation of the Quadratic Formula}
\author{Chenguang Zhang\thanks{Department of Mathematics, Massachusetts Institute of Technology. Email: \texttt{cgzhang@mit.edu}.}}

\date{January 5, 2020}

\begin{document}
\maketitle

\begin{abstract}
This article provides a simple geometric interpretation of the quadratic formula. The geometry helps to demystify the formula's complex appearance and casts it into a much simpler existence, thus potentially benefits early algebra students.
\end{abstract}

\section{Introduction}\label{introduction}

This work came from a note I wrote down several years back. Motivated by
the recent coverage of Po-Shen Lo's work~\cite{loh}, I decide it could
worth sharing. It is my hope that this short article will add yet another
pleasant aspect to the time-honored quadratic formula. To keep it short, the author will not review the rich history behind the quadratic formula, and instead highly recommend readers to consult Lo's article~\cite{loh} and references therein.

\section{Derivation}\label{derivation}

The quadratic formula gives the roots to the algebraic equation
$ax^{2} + bx + c = 0$. Without loss of generality, we set \(a = 1\) and assume \(b > 0\) and
\(c > 0\). So we are left with
\begin{equation}\label{eqn:quadratic}
x^{2} + bx + c = 0
\end{equation}
to solve. Let us first start with a simpler function (or curve)
\begin{equation}\label{eqn:curve}
y = x^{2} + bx.
\end{equation}

\begin{figure}
	\centering
	\includegraphics[width=0.75\linewidth]{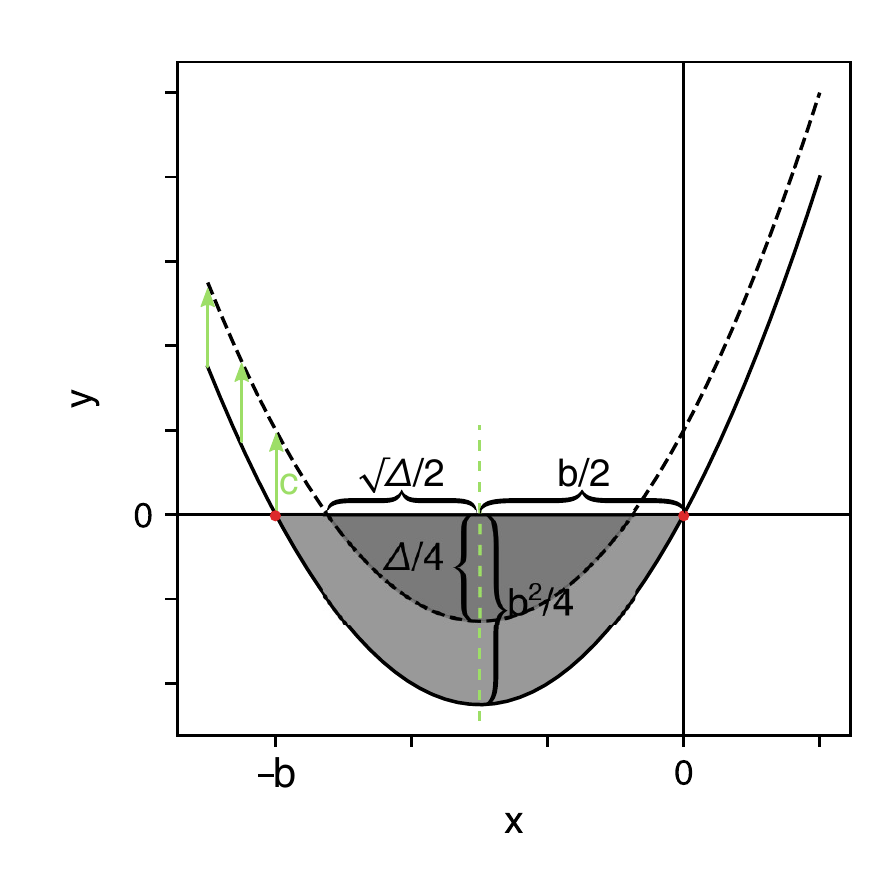}
	\caption{The curves \(y = x^{2} + bx\) (solid, associated
		with the light-shaded region) and \(y=x^{2} + bx + c\) (dashed,
		associated with the dark-shaded region). The green arrows indicate how the two curves are related: by a vertical shift.}
	\label{fig:illu}
\end{figure}

Algebraically, its roots are solutions to the equation \(y = x^{2} + bx = 0\), which are clearly \(x_{1} = 0\) and \(x_{2} = - b\). Geometrically, its roots
are the intersection of the curve with the \(x\)-axis, which are marked
by the two red dots in \cref{fig:illu}.

Before adding \(c\) back, we should realize two geometric properties
of curve \ref{eqn:curve} (see \cref{fig:illu}).
\begin{enumerate}
	\item Both the curve and its two roots are symmetric about \(x = - b/2\)
	(the vertical dashed line in green). This property is quite intuitive
	from \cref{fig:illu} and can be shown by averaging the two roots. When we add
	\(c\) back, we merely lift the curve vertically upward. Thus, the
	symmetry is respected.
	\item For easy reference, we refer to the region between the $x$-axis and the part of curve below it as the ``valley''. They are shaded in \cref{fig:illu}. The ``half-width''
        of this valley is $b/2$, and the ``depth'' of this valley is $b^2/4$. While the depth can be found by directly substituting \(x = - b/2\) into the equation, a better way is
        to use the quadratic nature of the curve: as one travels along the curve from its very bottom, the vertical shift equals the square of the horizontal shift \emph{by
        definition}. Naturally, the depth of the valley (here the vertical shift) equals the square of the \emph{half} width of the valley (here the horizontal shift). Conversely, the half-width of the valley equals the \emph{square root} of the valley depth.
\end{enumerate}

Now we add \(c\) back. By item 1, adding \(c\) lifts \(y = x^{2} + bx\)
vertically upward by an amount of \(c\), so the valley is shallower and its
new depth is \(\frac{b^{2}}{4} - c\). By item 2, we know that the half
width of the new valley is \(\sqrt{\frac{b^{2}}{4} - c}\). Because the two new roots are still symmetric about \(x = - b/2\), we have the roots
\begin{equation}\label{eqn:root1}
-\frac{b}{2} \pm \sqrt{\frac{b^{2}}{4} - c},
\end{equation}
or in the more familiar form:
\begin{equation}\label{eqn:root2}
\frac{- b \pm \sqrt{b^{2} - 4c}}{2}.
\end{equation}

\section{Discussion}\label{discussion}

\subsection{Alternative split}
In this article, the author splits \cref{eqn:quadratic} into $x^2+bx=0$ and $+c$ to derive the geometric interpretation. He also tried other ways of splitting. For example, $x^2+c$ and $bx$, or $x^2$ and $bx+c$. They failed to work due to that the term $bx$ or  $bx+c$ represents a slanted line, making it hard to proceed. It is likely that the split used by this article is the easiest one that yields to clear geometric interpretations.

\subsection{The discriminant}\label{the-Discriminant}
An interesting observation is that the discriminant of the quadratic
equation \(\Delta = b^{2} - 4c\) is nothing but the full width of
the valley \emph{squared}. Certainly, there will be no roots when this value is negative. Yet another interpretation is that the depth of the valley is
\(\frac{1}{4}\Delta=\frac{b^{2}}{4} - c\). When it is negative, the valley is
non-existent, the curve is above the $x$-axis (by an amount of $-\frac{1}{4}\Delta$), and there are no roots.

At this point, it is clear that a positive \(c\) always tries to eliminate the valley and roots. Whereas \(b\), regardless of its sign, always creates two roots. In this sense, the two coefficients compete against each other.

\subsection{Pedagogy}\label{pedagogy}
For the students to follow this article, it is vital that they first comprehend the geometry of quadratic curves. The part likely needs the most explanation is item 2 of the previous section, which relates the horizontal and vertical changes of the curve, as one travels along it from its very bottom.

\begin{figure}
	\centering
	\includegraphics[width=\linewidth]{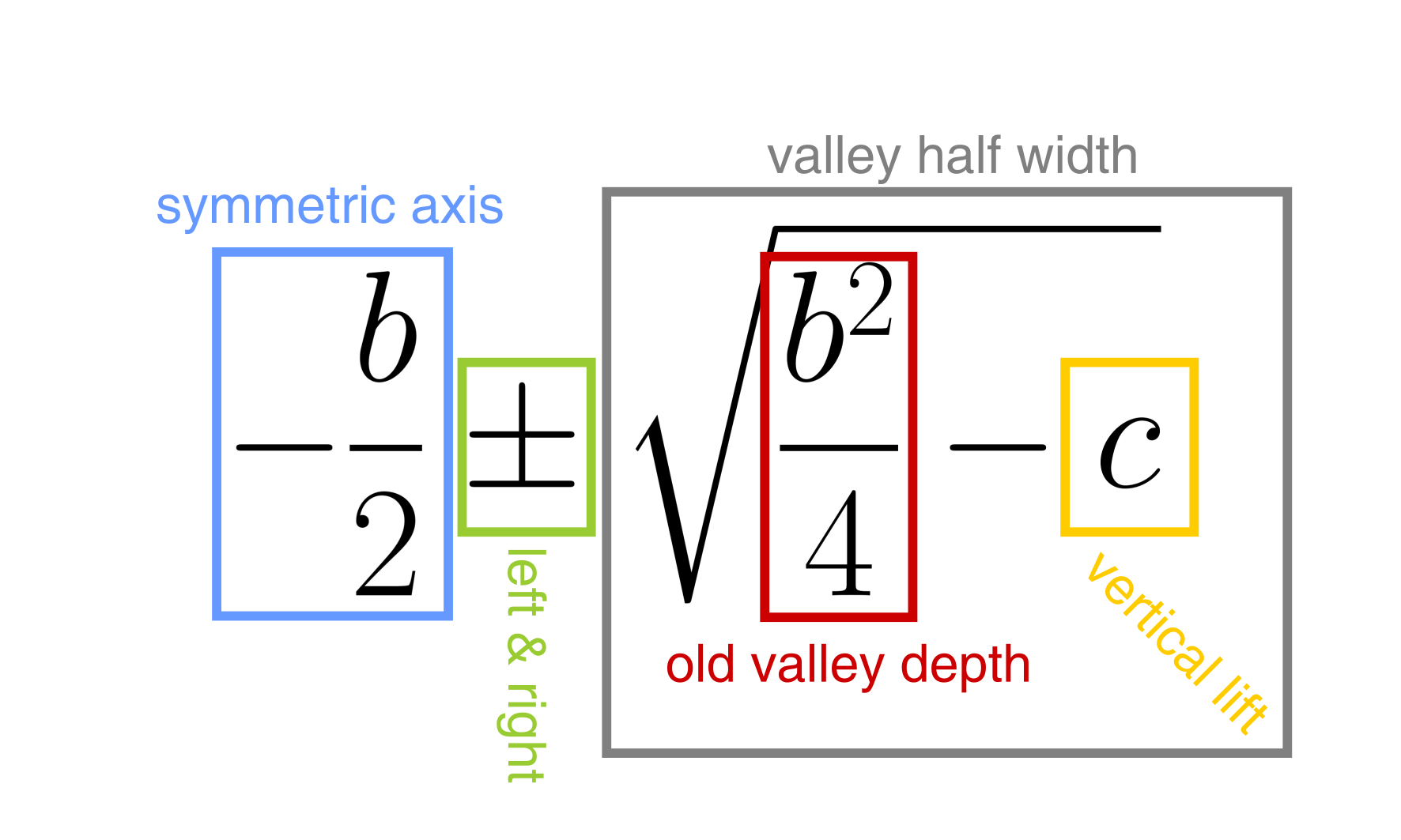}
	\caption{\cref{eqn:root1}, with each term annotated by its geometric meaning .}
	\label{fig:illu2}
\end{figure}

As a final comment. \cref{eqn:root1} seems to have clearer meaning than \cref{eqn:root2}. In the author's opinion, there are clear geometric interpretations of each value and operator in \cref{eqn:root1}; a heavily annotated version of which is shown in \cref{fig:illu2}. Much of these interpretations are lost in \cref{eqn:root2}, leaving a lifeless instruction waiting to be executed mechanically.

\end{document}